\documentclass{amsart}

\usepackage[german, english]{babel}
\usepackage{amssymb,amsmath,amsthm,amscd}

\advance\oddsidemargin by -0.1cm
\advance\evensidemargin by -0.7cm
\newtheorem{theorem}{Theorem}
\newtheorem{proposition}{Proposition}
\newtheorem{corollary}{Corollary}
\newtheorem{lemma}{Lemma}
\newtheorem{definition}{Definition}
\theoremstyle{definition}
\newtheorem{remark}{Remark}
\newtheorem{example}{Example}

\newcommand{\modulus}[1]{\left\lvert #1 \right\rvert}
\newcommand{\norm}[1]{\left\| #1 \right\|}
\newcommand{\mklm}[1]{\left\{ #1 \right\}}
\newcommand{\eklm}[1]{\left\langle #1 \right\rangle}

\renewcommand{\d}{\,d}

\newcommand{\C}{{\mathbb C}}
\newcommand{\R}{{\mathbb R}}
\newcommand{\A}{{\mathcal A}}
\newcommand{\B}{{\mathcal B}}
\newcommand{\D}{{\mathcal D}}

\newcommand{\I}{{\mathcal I}}

\renewcommand{\P}{{\mathcal P}}

\newcommand{\Af}{{\mathfrak A}}
\newcommand{\Bf}{{\mathfrak B}}

\newcommand{\1}{{\bf 1}}
\renewcommand{\epsilon}{\varepsilon}
\renewcommand{\phi}{\varphi}
\renewcommand{\rho}{\varrho}

\newcommand{\Cinft}{{\rm C^{\infty}}}
\newcommand{\Ctest}{{\rm C^{\infty}_c}}
\newcommand{\Cvan}{{\rm C_0}}

\renewcommand{\L}{{\rm L}}

\newcommand{\GL}{\mathrm{GL}}
\newcommand{\SL}{\mathrm{SL}}
\newcommand{\SO}{\mathrm{SO}}

\renewcommand{\b}{{\bf \mathfrak b}}
\newcommand{\g}{{\bf \mathfrak g}}
\newcommand{\h}{{\bf \mathfrak h}}
\renewcommand{\k}{{\bf \mathfrak k}}
\newcommand{\n}{{\bf \mathfrak n}}
\newcommand{\p}{{\bf \mathfrak p}}
\renewcommand{\t}{{\bf \mathfrak t}}
\newcommand{\U}{{\bf \mathfrak U}}

\newcommand{\Ad}{\text{Ad}\,}

\renewcommand{\det}{\text{det}\,}

\DeclareMathOperator{\gd}{\partial}

\DeclareMathOperator{\grad}{grad}

\DeclareMathOperator{\slim}{s-lim}
\newcommand{\e}[1]{\,{\mathrm e}^{#1}\,}

\begin{document}

\author{Pablo Ramacher}
\title{Analysis on real affine $G$--varieties}
\address{Pablo Ramacher, Humboldt--Universit\"at zu Berlin, Institut f\"ur Reine Mathematik, Sitz: Rudower Chaussee 25, D--10099 Berlin, Germany}
\email{ramacher@mathematik.hu-berlin.de}
\thanks{The author wishes to thank Prof. Thomas Friedrich for suggesting the initial line of research of this paper. This work was supported by the SFB 288 of the DFG}

\begin{abstract} We consider the action of a real linear algebraic group $G$ on a smooth, real affine algebraic variety $M\subset \R^n$, and study the corresponding left regular $G$-representation on the Banach space $\Cvan(M)$ of continuous, complex valued functions on $M$ vanishing at infinity. We show that the differential structure of this representation is already completely characterized by the action of the Lie algebra $\g$ of $G$ on the dense subspace $\P=\C[M] \cdot
e^{-r^2}$, where $\C[M]$ denotes the algebra of regular functions of $M$ and $r$ the distance function in $\R^n$. We prove that the elements of this subspace constitute analytic vectors of the considered $G$-representation, and, using this fact, we construct discrete reducing series in $\Cvan(M)$.  In case that $G$ is reductive, $K$ a maximal compact subgroup, $\P$ turns out to be a $(\g,K)$-module in the sense of Harish-Chandra and Lepowsky, and by taking suitable subquotients of $\P$, respectively $\Cvan(M)$,  one gets admissible  $(\g,K)$-modules as well as $K$-finite Banach representations. 
\end{abstract}

\maketitle

 \tableofcontents

\section{Introduction}

In the present article, we consider the regular action of a real linear algebraic group $G\subset \GL(n,\R)$  on a smooth, real affine algebraic variety $M\subset \R^n$, and the corresponding continuous, left regular  representation $\pi$ of $G$ on the Banach space $\Cvan(M)$ of continuous, complex valued functions on $M$ vanishing at infinity. Fix a compact subgroup $K$ of $G$. In general representation theory, a crucial role is played by the space of differentiable, $K$-finite vectors. If $E$ denotes a locally convex, complete, Hausdorff, topological vector space, and $\sigma$ a continuous representation on $E$, it is defined as the algebraic sum
\begin{displaymath}
  E_K=\sum _{\lambda \in \hat K} E_\infty \cap E(\lambda),
\end{displaymath}
where $E_\infty$ is the space of differentiable vectors in $E$ for $\sigma$, $\hat K$ the set of all equivalence classes of finite dimensional irreducible representations of $K$, and $E(\lambda)$ the isotypic $K$-submodule of $E$ of type $\lambda \in \hat K$. $E_K$ is dense in $E$, and the study of $\sigma$ can be reduced to a great extent to the study of the module $E_K$.
In case of the Banach representation $(\pi,\Cvan(M))$, a more natural submodule associated with the algebraic structure of $M$ arises. Let $\C[M]$ denote the ring of regular functions of $M$,  $r$ the distance function in $\R^n$, and consider the subspace
$$\P=\C[M] \cdot e^{-r^2}\subset \Cvan(M).$$
It was introduced by Agricola and Friedrich in \cite{agricola-friedrich}, and they proved that it is dense in $\Cvan(M)$. Let $\g$ denote the Lie algebra of $G$, and $\U(\g_\C)$ the universal enveloping algebra of the complexification of $\g$. In what follows, we will show that the differential structure of the representation $\pi$ is already completely determined by the action of $\g$ on $\P$. We prove that $\P$  is contained in the space of analytic elements of $\Cvan(M)$, which allows us to derive $\pi(G)$-invariant decompositions of the Banach space $\Cvan(M)$ from algebraic decompositions of $\P$ into $d\pi(\U(\g_\C))$-invariant subspaces. In particular, we obtain discrete reducing series in $\Cvan(M)$. In case that $G$ is reductive, $\P$ turns out to be a $(\g,K)$-module in the sense of Harish-Chandra and Lepowsky, and one has $\P=\sum_{\delta \in \hat K} \P \cap \Cvan(M)(\lambda)$. By taking suitable subquotients of $\P$ and $\Cvan(M)$, one gets admissible $(\g,K)$-modules as well as $K$-finite Banach representations.

\section{Some general remarks on  Banach representations} 
\label{sec:1}

Let us  begin with some generalities concerning Banach representations.
Thus, consider a weakly continuous representation $\pi$ of a Lie group $G$ on a Banach space 
$\B$ by linear bounded operators. According to a  result of Yosida from the theory of  $(C_0)$-semigroups \cite{yosida},  $\pi$ is also continuous with respect to the strong topology on $\B$, so that both continuity concepts coincide. Such a representation is called a \emph{Banach representation}. Similarly, the generators of the weakly, respectively      strongly,  continuous one-parameter groups of operators $h \mapsto \pi({\e{hX}})$ coincide, where
$h \in \R$, and $X$ is an element of the Lie-algebra $\g$
of   $G$. We will denote the corresponding generators by $d\pi(X)$, which are given explicitly by 
\begin{displaymath}
d\pi(X) \phi = \frac d {dh} \pi (\e{hX}) \phi_{\big|_{h=0}}
\end{displaymath}
for those $\phi \in \B$, for which the limit exists. These operators are closed and densely defined with respect to the weak and strong topology on $\B$. Let $d_G$ be left invariant Haar-measure on $G$,  $f \in
\L^1(G,d_G)$, and $\B^\ast$ the dual of $\B$.  Then, for each $\phi \in \B$, $f(g)
\mu(\pi(g)\phi)$ is $d_G$-integrable for arbitrary $\mu \in \B^\ast$,
and there exists a $\psi \in \B$ such that  
\begin{displaymath}
\mu(\psi)=\int _G f(g) \mu(\pi(g) \phi) \d_G(g)
\end{displaymath}
for all $\mu \in \B^\ast$, i.\ e.\ $f(g) \pi(g) \phi$ is integrable in the sense of Pettis, see e.\ g.\
\cite{hille-phillips}. Here $\psi$ is given as a weak limit, and one defines as this limit the integral
\begin{displaymath}
\int_G f(g) \pi(g) \phi \d_G(g)=\psi,
\end{displaymath}
in this way getting a linear operator $\pi(f): \B \rightarrow
\B, \phi \mapsto \int_G f(g) \pi(g) \phi \d_G(g)$. Note that
    $\norm{\pi(f)}\leq \norm{f}_{\L^1}$. In case that $f$ is $\L^1$-integrable and continuous, $f(g) \pi(g) \phi$ is also integrable in the sense of Bochner,  and $\psi$ given directly by the corresponding strongly convergent integral. Let
\begin{displaymath}
\B_\infty=\mklm{\sum_{i=1}^l \pi(f_i) \phi_i: \quad f_i \in
\Ctest(G),\, \phi_i \in \B,\, l=1,2,3,\dots}  
\end{displaymath}
be the  G{\aa}rding-subspace of $\B$ with respect to  $\pi$, and $\B_{\infty,s}$, respectively\ $\B_{\infty,w}$, the subspace of differentiable elements in $\B$ with respect to the strong, respectively\ weak, topology. $\B_\infty$ is norm-dense 
in $\B$ and,  according to 
Langlands \cite{langlands}, the generators $d\pi(X)$ are already completely determined by their action on the 
G{\aa}rding-subspace. Thus, if 
$\Gamma_{X_1,\dots,X_k}$ denotes the graph of the generators $d\pi(X_i)$,
$i=1,\dots,k$, and $\Gamma_{X_1,\dots,X_k|\B_\infty}$ its restriction to $\B_\infty$, one has
$\Gamma_{X_1,\dots,X_k}=\overline{\Gamma_{X_1,\dots,X_k|\B_\infty}}$.
As an immediate consequence, the differential structure of the representation $\pi$ is completely characterized by the action of the operators $d\pi(X)$ on $\B_\infty$. In particular, this implies that the strongly, respectively\ weakly, differentiable elements in $\B$ do coincide with those that are differentiable with respect to the one-parameter groups of operators $h \mapsto \pi({\e{hX}})$, the underlying topology being the strong, respectively\ weak, topology. Since, by Yosida, strong and weak generators
coincide, one finally has
\begin{equation}
\label{eq:0}
\B_{\infty,s}=\bigcap_{i=1}^d  \bigcap_{k\geq1} \D(d\pi(a_i)^k)= \B_{\infty,w},
\end{equation}
where $a_1,\dots a_d$ denotes a basis of 
$\g$, compare e.\ g.\ \cite{robinson}. On the other hand, by Dixmier and Malliavin \cite{dixmier-malliavin}, $\B_\infty$ coincides with $\B_{\infty,s}$, and therefore it follows that
\begin{displaymath}
\B_\infty=\B_{\infty,s}=\B_{\infty,w}.
\end{displaymath}
$\B_\infty$ is invariant under the $G$-action  $\pi$ and the $\g$-action $d\pi$, by which one gets representations of
$G$ and $\g$ on  $\B_\infty$. The fact that $\B_{\infty,s}=\B_{\infty,w}$ can also be deduced from the following general argument going back to Grothendieck: If $M$ is a non-compact $\Cinft$-manifold and $E$ a locally convex, complete, Hausdorff, topological vector space, then  $f:M \to E$ is a $\Cinft$-mapping with respect to the locally convex topology if and only if, for all $\mu\in E^\ast$, the function $\mu(f(m))$  on  $M$ is infinitely often differentiable, see \cite{warner}, page 484.

Assume now that $G$ is a real linear algebraic group.
As a smooth, real, affine algebraic variety, $G$ is a real analytic manifold,
and hence a real analytic Lie group. Therefore, the exponential map is, locally, a real analytic homeomorphism.
Taking a sufficiently small neighbourhood of zero in $\g$, and 
assuming a  decomposition of $\g$ of the form $\g_1 \oplus \dots \oplus \g_l$,
$(X_1, \dots , X_l)\mapsto g \e{X_1} \dots   \e{X_l}$ becomes  an analytic homeomorphism of the aforementioned neighbourhood onto an open neighbourhood of  $g\in G$. With  the identification $\g\simeq \R^d$, and with respect to a basis $a_1, \dots, a_d$ of $\g$, the canonical coordinates of second type of a point $g\in G$ are then given by 
\begin{equation}
\label{eq:0a}
\Phi_g:g U_e\ni g \e{t_1a_1}\dots \e{t_da_d} \mapsto (t_1,\dots, t_d) \in W_0,
\end{equation}
where     $W_0$ denotes a sufficiently small neighbourhood of $0$ in $\R^d$, and $U_e=\exp(W_0)$. We will write for $\Phi_e$ simply  $\Phi$. Let $\phi \in \B_\infty$, so that $g \mapsto \pi(g) \phi$
becomes a  $\Cinft$-map from $G$ to   $\B$ with respect to the strong and weak topology of $\B$. This is  equivalent to the fact that, for all $g \in G$, the 
map $(t_1,\dots,t_d) \mapsto \pi({\Phi_g^{-1}(t)}) \phi$ is infinitely often strongly, respectively weakly, differentiable on  $W_0$. With regard to any of these topologies, we obtain for $t  \in W_0$ the relations  
\begin{align}
\label{eq:1}
\begin{split}
d\pi(a_j) &\pi(\Phi_g^{-1}(t))\phi=\lim_{h \to 0} h^{-1}
[\pi(\e{ha_j})-\1]\pi({\Phi_g^{-1}(t)})\phi\\
&=\frac d {dh} \Big (\pi\big({\Phi_g^{-1}(s^j(h,t))\big)}\phi\Big )_{\big\vert_{h=0}}=\sum_{k=1}^d \frac \gd {\gd t_k} \Big (\pi\big(\Phi_g^{-1}(t)\big )
\phi\Big ) \frac
d {dh} s^j_k(0,t);
\end{split} 
\end{align}
here the $s_k^j(h,t)$ are real analytic functions such that 
$$
\e{ha_j} g \e{t_1a_1} \dots \e{t_da_d}=g \e{s_1^j(h,t)a_1}\dots \e{s_d^j(h,t)a_d},
$$
since $\e{ha_j} g \e{t_1a_1} \dots \e{t_da_d} \in \Phi_g^{-1}(W_0)$ for small $h$. In a similar way, we have
\begin{align}
\label{eq:2}
\pi(\Phi_g^{-1}(t))d\pi(a_j) \phi&=\sum_{k=1}^d \frac \gd {\gd t_k} \Big (\pi\big(\Phi_g^{-1}(t)\big )
\phi\Big ) \frac
d {dh} r^j_k(0,t),
\end{align}
with real analytic functions $r_k^j(h,t)$ satisfying the relations 
$$
g \e{t_1a_1} \dots \e{t_da_d}\e{ha_j} =g \e{r_1^j(h,t)a_1}\dots \e{r_d^j(h,t)a_d}.
$$
Clearly,  $s^j_k(0,t)=r^j_k(0,t)=t_k$, and
$s^j_k(h,0)=r^j_k(h,0)=\delta_{kj}h$ for  $g=e$. Finally, we also note that 
\begin{align}
\label{eq:3}
\begin{split}
\frac \gd {\gd t_j} \Big (\pi(\Phi_g^{-1}(t)) \phi\Big )&=\lim_{h \to 0} h^{-1}
\pi(\Phi_g^{-1}(t)) [\pi\big(\Phi_e^{-1}(t^j(h,t))\big)-\1]\phi\\
&=\sum_{k=1}^d \pi(\Phi_g^{-1}(t)) \d \pi(a_k) \phi \,\,\frac d {dh} t^j_k(0,t),
\end{split}
\end{align}
where the $t^j_k(h,t)$ are real analytic functions in $h$ and $t$, and satisfy the relations $$\e{-t_da_d} \dots
\e{-t_{j+1}a_{j+1}}\e{ha_j}\e{t_{j+1}a_{j+1}}\dots
\e{t_da_d}=\e{t_1^j(h,t)a_1}\dots \e{t_d^j(h,t)a_d}.$$ 
One has
$t^j_k(0,t)=0$, $t^j_k(h,0)=\delta_{kj}h$, so that, in particular,
$$
\frac \gd {\gd t_k} \left(\pi(\Phi_g^{-1}(t)) \phi\right )_{\big\vert_{t=0}}=\pi(g)d\pi(a_k)\phi.
$$

\section{The regular representation $(\pi,\Cvan(M))$ of a real algebraic group $G$ acting on a smooth, real affine variety $M$}
\label{sec:3}

We come now to the proper subject of this paper.
Let $M$ be a smooth, real affine algebraic variety and $G$ a real linear algebraic group, which acts regularly on $M$. In what follows, we will view $G$ as a closed subgroup of $\GL(n,\R)$, and $M$ as embedded in $\R^n$. Denote by $\Cvan(M)$ the vector space of all continuous, complex valued functions on $M$ which vanish at infinity; provided  with the supremum norm, $\Cvan(M)$ becomes a Banach space. According to the Riesz representation theorem for locally compact, Hausdorff, topological vector  spaces, its dual is given by the Banach space of all regular, complex measures $\mu:{\mathfrak{B}} \rightarrow \C$ on $M$ with norm $\vert \mu\vert (M)$. Here ${\mathfrak{B}}$ denotes the $\sigma$--algebra of all Borel sets of $M$, $\vert \mu\vert$ the variation of  $\mu$, and one has $\Cvan(M)\subset{\mathcal L}^1(\mu)$. The $G$-action on $M$ induces a representation $\pi$ of $G$ on $\Cvan(M)$ by bounded linear operators according to
\begin{displaymath}
\pi(g):\Cvan(M) \rightarrow \Cvan(M), \quad (\pi(g) 
\phi)(m)=\phi(g^{-1} m),
\end{displaymath}
where $g \in G$. Henceforth, this representation will be called the \emph{left regular representation} of  $G$ on  $\Cvan(M)$. It is continuous with respect to the weak topology on $\Cvan(M)$, which is characterized by the family of seminorms $\vert
\phi\vert_{\mu_1,\dots,\mu_l}=\sup_i \vert \mu_i(\phi)\vert$, $\mu_i
\in \Cvan(M)^\ast, \, l\geq 1$. Indeed, by the theorem of Lebesgue on bounded convergence, one immediately deduces  
\begin{displaymath}
\lim_{g \to e} \mu(\pi(g) \phi)=\lim_{g \to e} \int \phi(g^{-1}m) \d\mu(m)= \int \phi(m) \d\mu(m)=\mu(\phi)
\end{displaymath}
for all $\phi \in \Cvan(M)$ and complex measures $\mu$, as well as 
$\mu\circ \pi(g)\in \Cvan(M)^\ast$ for all $g \in G$. Hence, by the considerations of the previous section, $\pi$ is a Banach representation of $G$.

\subsection{The Gaussian measure on  affine varieties and the $\g$-module $\P=\C[M]\cdot e^{-r^2}$}

In the following, let    $\C[M]$ be the ring of all functions that arise by restriction  of polynomials in $\R^n$ to $M$ and denote by $\P$ the subspace 
\begin{displaymath}
\P=\C[M] \cdot e^{-r^2},
\end{displaymath}
where $r^2(m)=m_1^2+\dots+m_n^2=m^2$ is the square of  the distance of a point $m \in M$ to the origin in $\R^n$ with respect to the coordinates $m_1,\dots,m_n$.
According to  Agricola and Friedrich \cite{agricola-friedrich}, $\P$ is norm-dense in  $\Cvan(M)$. Although $\P$ is not invariant under the $G$-action $\pi$, the next proposition shows that $\P$ is a $d\pi(\U(\g_\C))$-invariant subspace of $\Cvan(M)_\infty$. As will be shown later, the elements in $\P$ are even analytic.
\begin{proposition}
\label{prop:1}
$\P$ is  a   $\g$-submodule of  $\Cvan(M)_\infty$.
\end{proposition}
\begin{proof}
As already explained, we can view $\Cvan(M)$ as endowed with the weak topology. Let $a_1, \dots, a_d$ be a basis of $\g$,  $p \in \C[M]$ a polynomial on $M$, and  
$\phi=p \cdot e^{-r^2}$ an element of  $\P\subset\Cvan(M)$. One computes
\begin{equation}
\label{eq:6c}
\frac d {dh} \Big (\phi(\e{-ha_j}m)\Big )_{\big \vert_{h=0}}=\eklm {(\grad \phi)(m),\frac d {dh} (\e{-ha_j}m)_{\big \vert_{h=0}}}, 
\end{equation}
as well as
\begin{displaymath}
(\grad \phi)(m)=\big ((\grad p)(m)-2 p(m) m\big) e^{-m^2}.
\end{displaymath}
By assumption,  $G\subset \GL(n,\R)$ and $\g \subset M_n(\R)$ act on
$M\subset \R^n$ by matrices, and we obtain 
\begin{displaymath}
\frac d {dh} \Big(\e{-ha_j}m\Big)_{\big \vert_{h=0}}=\lim_{h \to 0} h^{-1}
\sum\limits_{k=1}^{\infty} \frac{(-h a_j)^k} {k!}m=-a_j m.
\end{displaymath}
Hence, if   $\tilde a_j$ denotes the vector field $(\tilde
a_j)_m=a_j m$ on $M$, one has $-\tilde a_j \phi=d/dh\,\phi(\e{-ha_j}\cdot)_{\vert h=0}$ $\in
\P\subset \Cvan(M)$.
For arbitrary $\mu\in \Cvan(M)^\ast$, we therefore get, according to Lebesgue,
\begin{align*}
\lim_{h \to 0} h^{-1} \mu \big (
[\pi(\e{ha_j})-\1]\phi\big)&=\int_M\lim_{h \to 0} h^{-1} \big [ \phi(\e{-ha_j} m)
-\phi(m)\big ] \d\mu(m)
=-\mu (\tilde a_j \phi ).
\end{align*}
Hence $\phi \in \D(d\pi(a_j))$ for all  $j=1, \dots, d$, and $d\pi(a_j)
\phi=-\tilde a_j \phi$. Since $\tilde a_j \phi \in \P$, the assertion now follows with \eqref{eq:0}.
\end{proof}
As a consequence, the relations \eqref{eq:1}--\eqref{eq:3}
hold for $\phi \in \P$, also.

\subsection{An approximation argument} 
\label{sec:2.2}

Let $a_1,\dots,a_d$ be a basis of $\g$ as above. In the following we will denote the generators
$d\pi(a_i)$ by  $A_i$. We set 
\begin{displaymath}
\Gamma_{a_1,\dots,a_d|\P}=\mklm{(\phi,A_1\phi,\dots,A_d \phi) \in \Cvan(M) \times\dots \times
\Cvan(M):\phi \in \P}.
\end{displaymath}
As already noted, $\P$
is not $\pi(G)$-invariant in general, so that, for  $\phi \in \P$, it is not true that 
$(\pi(g) \phi,A_1\pi(g)\phi,\dots,A_d \pi(g)\phi)$ is contained in $\Gamma_{a_1,\dots,a_d|\P}$. Nevertheless, it will be shown in the following that
the last assertion is  correct, if instead of $\Gamma_{a_1,\dots,a_d|\P}$ one considers its closure $\overline{ \Gamma_{a_1,\dots,a_d|\P}}$ in $\Cvan(M)
\times\dots \times
\Cvan(M) $ with respect to the strong product topology, and $g$ is taken in a sufficiently small neighbourhood of the unit. To start with, we will need the following lemma.
\begin{lemma}
\label{lemma:1}
Let $\rho>0$ and $\kappa \in \R^\ast$. Then, for   $l \to
\infty$,
\begin{displaymath}
e^{-\rho x^2} \sum_{k=0}^{l-1} \frac{(\kappa x^2)^k}{k!}\longrightarrow
e^{-\rho x^2} e^{\kappa x^2}
\end{displaymath}
converges uniformly on $\R^n$, provided $\rho/\vert \kappa \vert \geq2$.
\end{lemma}
\begin{proof}
By the comparison criterium for series, one has 
\begin{displaymath}
\Big |e^{\kappa x^2}-\sum_{k=0}^{l-1} \frac{(\kappa
x^2)^k}{k!}\Big |\leq \frac {(|\kappa | x^2)^l}{l!} e^{|
\kappa | x^2},
\end{displaymath}
thus obtaining
\begin{displaymath}
\sup_{x \in \R^n} \Big | e^{-\rho x^2}\Big (e^{\kappa x^2}-\sum_{k=0}^{l-1} \frac{(\kappa
x^2)^k}{k!}\Big ) \Big |\leq \sup_{x \in \R^n} \frac {(|\kappa | x^2)^l}{l!} e^{x^2(|
\kappa |-\rho)}.
\end{displaymath}
Assume $\rho/|\kappa|>1$, so that $\lambda=\rho/|\kappa|-1>0$. The supremum on the right hand side of the last inequality, denoted in the following by $\Sigma_{l,\lambda}$, can be computed as
\begin{equation}
\label{eq:6b}
\Sigma_{l,\lambda}=\frac 1 {l!} \sup_{y\geq 0}\,  {y^l}e^{-\lambda y}=\frac 1
{l!}  (l/\lambda)^l \, e^{-l}.
\end{equation}
By the Stirling formula, $l!$ is asymptotically given by  $\sqrt{l}\,\,l^l
e^{-l}$, so that, in case $\lambda \geq 1$,  
one deduces    $\Sigma_{l,\lambda} \to 0$ for  $l \to \infty$, and hence the assertion.
\end{proof}
Now, we are able to proof the announced result. For this sake, let us define the set
\begin{displaymath}
\Gamma_{a_1,\dots,a_d|\pi(U) \P}=\mklm{(\phi, A_1\phi,\dots, A_d\phi) \in \Cvan(M) \times \dots \times \Cvan(M): \phi \in \pi(g)\P,\, g \in U},
\end{displaymath}
where $U$ denotes a neighbourhood of $e\in G$.
\begin{proposition}
\label{prop:2}
Let $U$ be a sufficiently small neighbourhood of the unit in $G$, $\phi \in \P$, $g\in U$. Then there exists a series $\phi^g_k \in \P$ such that $\norm{\pi(g)\phi-\phi^g_k} \to 0$, and $\norm{A_j\pi(g)\phi-A_j\phi^g_k} \to 0$ for all $j=1,\dots,d$. In other words, 
\begin{displaymath}
\Gamma_{a_1,\dots,a_d|\pi(U)\P} \subset \overline {\Gamma_{a_1,\dots,a_d|\P}}.
\end{displaymath}
\end{proposition}
\begin{proof}
Let $\phi=p\cdot e^{-r^2}\in \P$, and define on the subspace $\P$ the linear operators
\begin{equation}
\label{eq:6a}
\pi_k(g) \phi=\rho(g)p\sum_{j=0}^k \frac {(r^2-(\rho(g)r^2)^j}{j!} \cdot e^{-r^2}, \quad g \in G,
\end{equation}
where $\rho(g)p(m)=p(g^{-1}m)$ denotes the left regular representation of $G$ on $\C[M]$. Note that $\P$ is left invariant under the operators \eqref{eq:6a}, which, in general, can not be extended to continuous operators on   $\Cvan(M)$. We set   $\phi^g_k=\pi_k(g)\phi$, and compute 
\begin{align*}
\|\pi(&g)\phi-\phi^g_{k-1}\|=\sup_{m\in M}\Big \vert (\rho(g) p)(m) \Big [ e^{-(g^{-1}m)^2}-e^{-m^2}\sum_{j=0}^{k-1}\frac 1 {j!} (m^2-(g^{-1}m)^2)^j  \Big ] \Big \vert\\
&\leq \sup_{m\in M} \Big \vert (\rho(g) p)(m) e ^{-m^2/2}\Big \vert\sup_{m\in M} \Big \vert e^{-m^2/2}\Big ( e^{m^2-(g^{-1}m)^2}-\sum_{j=0}^{k-1} \frac 1 {j!} (m^2-(g^{-1}m)^2)^j\Big ) \Big\vert.
\end{align*}
Writing $g_{ij}$ for the matrix entries of an arbitrary element $g\in G$, one computes
\begin{align*}
(gm)^2&=(g_{11}m_1+\dots + g_{1n}m_n)^2+\dots +(g_{n1}m_1+\dots + g_{nn}m_n)^2\\
&=\sum^n_{k,i=1} g_{ki}^2\,m_i^2+2 \sum ^n_{k=1}\sum^n_{i<j} g_{ki} \,g_{kj} \,m_im_j\\
&\leq \sum_{i=1}^n g_{ii}^2 \, m_i^2+ \sum_{i\not=k}^n g_{ki}^2 m_i^2+\sum_{k=1}^n \sum_{i<j}^n g_{ki} \, g_{kj} (m_i^2+m_j^2)\\
&\leq \big ( \max_i g_{ii}^2+\max _{k\not=i} g_{ki}^2+\max_{k;i<j}\vert g_{ki} \, g_{kj} \vert\big) \, m^2;
\end{align*}
in particular,
\begin{equation*}
\vert (gm)^2-m^2 \vert \leq \big ( \max_i \vert g_{ii}^2-1\vert +\max _{k\not=i} g_{ki}^2+\max_{k;i<j}\vert g_{ki} \, g_{kj} \vert\big) \, m^2.
\end{equation*}
Setting $\kappa_{g^{-1}}=\max_i \vert g_{ii}^2-1\vert +\max _{k\not=i} g_{ki}^2+\max_{k;i<j}\vert g_{ki} \, g_{kj} \vert$, we therefore obtain the estimate 
\begin{equation}
\label{eq:6d}
 \Big\vert e^{m^2-(g^{-1}m)^2}-\sum_{j=0}^{k-1} \frac 1 {j!} (m^2-(g^{-1}m)^2)^j\Big \vert\leq \sum_{j\geq k}^\infty \frac 1 {j!} \vert m^2-(g^{-1}m)^2 \vert^j\leq \sum_{j\geq k}^\infty \frac 1 {j!} (\kappa_g m^2)^j,
\end{equation}
so that, with $C=\sup_{m\in M} \Big \vert (\rho(g) p)(m) e ^{-m^2/2}\Big \vert$,
\begin{equation*}
\|\pi(g)\phi-\phi^g_{k-1}\|\leq  C \sup_{m\in M} \Big \vert e^{-m^2/2}\Big ( e^{\kappa_g m^2}-\sum_{j=0}^{k-1} \frac 1 {j!} (\kappa_gm^2)^j\Big ) \Big\vert.
\end{equation*}
Assume that  $U$ is a sufficiently small neighbourhood of the unit in $G$ such that $\kappa_g \leq 1/4$ for all   $g\in U$. Then, by Lemma     \ref{lemma:1}, $\|\pi(g)\phi-\phi^g_{k-1}\|\rightarrow 0$ goes to zero 
 as  $ k \to \infty$, for arbitrary $g \in U$. It remains to show that,  for $j=1,\dots,d$, $\|A_j\pi(g)\phi-A_j\phi^g_{k-1}\|$ goes to zero as  $k$ goes to infinity. Now, since $A_j \phi=-\tilde a_j \phi=d/\d h\, \phi(\e{-ha_j}\cdot)_{|_{h=0}}$, one has
\begin{gather*}
\|A_j\pi(g)\phi-A_j\phi^g_{k-1}\|=\sup_M\Big\vert \frac d{dh} (\pi(g)\phi)(\e{-ha_j}m)_{|h=0}-\frac d{dh}(\pi_{k-1}(g) \phi)(\e{-ha_j}m)_{\big \vert_{h=0}}\Big \vert\\
\leq \sup_{m \in M} \Big \vert \frac d {dh} (\pi(g)\phi)(\e{-ha_j}m) _{\big \vert_{h=0}}\,e^{(g^{-1}m)^2-m^2}\Big (e^{m^2-(g^{-1}m)^2} -\sum_{j=0}^{k-1} \frac 1 {j!} (m^2-(g^{-1}m)^2)^j\Big ) \Big |\\
+\sup_{m \in M} \Big \vert (\pi(g)\phi)(m) \frac d {dh} \Big (e^{-q(h)}\sum_{j=0}^{k-1} \frac 1 {j!} (q(h))^j\Big )_{\big \vert_{h=0}} \Big |, 
\end{gather*}
where we set $q(h)=(\e{-ha_j}m)^2-(g^{-1}\e{-ha_j}m)^2$.
The first summand on the right hand side converges to zero for $k \to \infty$ since, as a consequence of \eqref{eq:6d}, and 
\begin{align*}
\frac d {dh} (\pi(g)\phi)(\e{-ha_j}m)_{\big |_{ h=0}}&= -\eklm {(\grad p)(g^{-1}m)-2\,
p(g^{-1}m) \, g^{-1}m, g^{-1}a_jm} e^{-(g^{-1}m)^2},
\end{align*}
it can be estimated by
\begin{displaymath}
C' \sup_M \Big \vert e^{-m^2/2}\Big ( e^{\kappa_g m^2}-\sum_{j=0}^{k-1} \frac 1 {j!} (\kappa_gm^2)^j\Big ) \Big\vert
\end{displaymath}
with  $C'=\sup_{m \in M} \Big \vert \eklm {(\grad
p)(g^{-1}m)-2\,p(g^{-1}m) \, g^{-1}m, g^{-1}a_jm}e^{-m^2/2}\Big\vert,$ and a repea\-ted application of Lemma \ref{lemma:1} yields the assertion. In order to estimate the second summand, we note that
\begin{align*}
 \frac d {dh} \Big (e^{-q(h)}\sum_{j=0}^{k-1} \frac 1 {j!}
 (q(h))^j\Big )_{|h=0}&= \frac d {dq} \Big (e^{-q}\sum_{j=0}^{k-1}
 \frac 1 {j!} q^j\Big )_{|q=m^2-(g^{-1}m)^2}\, \dot q(0)\\
&=-e^{(g^{-1}m)^2-m^2} \frac 1 {(k-1)!} (m^2-(g^{-1}m)^2)^{k-1} \,
 \dot q(0).
\end{align*}
Here $\dot q (0)$ denotes the polynomial in $m$ which is explicitly given by $
\dot q(0)=-2 (\eklm{m, a_jm}-\eklm{g^{-1}m,g^{-1}a_jm})$. We get for the second summand the upper bound
\begin{displaymath}
C'' \frac 1 {(k-1)!}\sup_{m \in M} e^{-m^2/2}  (\kappa_gm^2)^{k-1}=C''\Sigma_{k-1,\lambda},
\end{displaymath}
where $C''=\sup_{m\in M}\vert(\rho(g)p)(m)\dot q(0) e^{-m^2/2}\vert$ , and  $\lambda=1/2\kappa_g$; 
$\Sigma_{k-1,\lambda}$ was 
defined in \eqref{eq:6b}. By repeating the arguments given in the proof of Lemma \ref{lemma:1}, it follows that, for $g \in
U$, $U$ as above, the second summand also converges to zero for $k\to\infty$. We get
\begin{displaymath}
\norm{A_j(\pi(g)\phi-\phi^g_{k-1})} \rightarrow 0
\end{displaymath}
for   $k \to \infty$, $g \in U$, and   $j=1,\dots,d$. This proves the proposition.
\end{proof}
\begin{remark}
Proposition \ref{prop:2} implies  that, for $g\in G$ in a sufficiently small neighbourhood   $U$ of  $e$,
$ \pi(g) \phi$ has an expansion as an  absolutely convergent series 
in $\Cvan(M)$ given by
\begin{equation}
\label{eq:6}
 \pi(g) \phi=  \sum_{l=0}^\infty
\rho(g)p \,(r^2-\rho(g)r^2)^l/l! \cdot e^{-r^2},
\end{equation}
where $\rho$ is the left regular representation of $G$ on $\C[M]$. Indeed, it follows from the proof of Lemma \ref{lemma:1} that 
\begin{align*}
\norm{\rho(g)p \frac {(r^2-\rho(g)r^2)^l}{l!} e^{-r^2}} &\leq \frac {\kappa_g^l} {l!}  \sup_M |p(g^{-1}(m)) e^{-m^2/2}| \cdot \sup_M |m^{2l} e^{-m^2/2}|\\&=\frac {\kappa_g^l} {l!}  C \cdot \sup_{y\geq 0} e^{-y/2} y^l=\frac {\kappa_g^l} {l!}  C (2l)^l  e^{-l}\approx (2\kappa_g)^l C \frac 1 {\sqrt l},
\end{align*}
where we set $C=\sup_M |p(g^{-1}(m)) e^{-m^2/2}|$. For the definition of $\kappa_g$, see the proof of Proposition \ref{prop:2}.
\end{remark}

\subsection{A dense graph theorem for the generators $d\pi(X)$}

We are now in position to show that the generators $d\pi(X)$ are already completely determined by their restriction to the subspace $\P=\C[M] \cdot e^{-r^2}$. A similar  dense graph theorem  involving the G{\aa}rding-subspace was conjectured originally by Hille within the theory of strongly continuous semigroups \cite{hille-phillips}, and afterwards proved by Langlands in his doctoral thesis \cite{langlands}. Our proof follows essentially the one of Langlands, which, nevertheless, makes use of the $G$-invariance of the  G{\aa}rding-subspace. The fact that $\P$ is not $\pi(G)$-invariant will be overcome by means of theapproximation argument exhibited  in the previous section.
\begin{theorem}
Let $X_1,\dots, X_k\in \g$, and denote by $\Gamma_{X_1,\dots,X_k}$ the graph of the generators $d\pi(X_1),\dots,d\pi(X_k)$, i.\ e.\ the set
\begin{displaymath}
\label{thm:1}
\Big\{(\phi,d\pi(X_1)\phi, \dots, d\pi(X_k)\phi) \in
\Cvan(M) \times \dots \times\Cvan(M): \phi \in \bigcap_{i=1}^k \D(d\pi(X_i))\Big\}.
\end{displaymath}
Write ${\Gamma_{X_1,\dots,X_k}}_{|\P}$ for its restriction to $\P\times \dots \times \P$. Then, with respect to the strong product topology on    $\Cvan(M)\times \dots \times \Cvan(M)$, 
\begin{displaymath}
\Gamma_{X_1,\dots,X_k}=\overline{{\Gamma_{X_1,\dots,X_k}}_{|\P}}.
\end{displaymath}
In particular, $\Gamma_{a_1,\dots,a_d}=\overline{\Gamma_{a_1,\dots,a_d|\P}}$.
\end{theorem}
\begin{proof}
We assume in the following that $\mklm{a_q,\dots,a_d}$ represents a maximal linear independent subset of $\mklm{X_1,\dots,X_k}$. It suffices to verify the assertion for this set, then. Let 
$(\phi, \psi_q,\dots,\psi_d)\in
\overline{{\Gamma_{a_q,\dots,a_d}}_{|\P}}$. There exists a series of functions $\phi_n\in \P$ such that   
\begin{displaymath}
(\phi_n,A_q\phi_n,\dots,A_d\phi_n) \rightarrow (\phi, \psi_q,\dots,\psi_d)
\end{displaymath}
with respect to the strong topology in  $\Cvan(M)\times\dots \times \Cvan(M)$. We obtain $\psi_i=A_i\phi$ for all 
$i=q,\dots,d$,  the $A_i$ being norm-closed, and, thus, $(\phi, \psi_q,\dots,\psi_d)\in \Gamma_{a_q,\dots,a_d}$.
To prove the converse inclusion $\Gamma_{a_q,\dots,a_d}\subset
\overline{{\Gamma_{a_q,\dots,a_d}}_{|\P}}$, we consider local regularizations of $\pi$. Consider the canonical coordinates of second type $\Phi:U_e\rightarrow W_0$ introduced in \eqref{eq:0a}, and let $Q_\epsilon=[0,\epsilon]\times\dots \times [0,\epsilon]$ be a cube in $W_0$ of length $\epsilon$. If $dt$ stands for Lebesgue-measure on $\R^d$ and 
$ \chi_{Q_\epsilon}$ for the characteristic function of $Q_\epsilon$, then, by the weak continuity of the $G$-representation $\pi$, the map     $\chi_{Q_\epsilon}(t)\,\pi({\Phi^{-1}(t)})\phi $ is weakly measurable with respect to
$dt$ for $\phi \in \Cvan(M)$, as well as separable-valued, and therefore, by Pettis, strongly measurable, see e.\ g.\ \cite{hille-phillips}. Because of
\begin{displaymath}
\int_{Q_\epsilon} \norm{\pi({\Phi^{-1}(t)})\phi} \d t = \norm{\phi} \cdot \epsilon^d,
\end{displaymath}
it follows that $\chi_{Q_\epsilon}(t)\,\pi({\Phi^{-1}(t)})\phi $ is
Bochner-integrable, and we define on $\Cvan(M)$ the linear bounded operators
$\pi(\chi_{Q_\epsilon})\phi=\int _{Q_\epsilon} {\pi({\Phi^{-1}(t)})\phi} \d t
$, in accordance with the regularizations $\pi(f)$ of  $\pi$ already introduced. Clearly, $\norm{\epsilon^{-d} 
\pi(\chi_{Q_\epsilon})\phi-\phi}\to 0$ for  $\epsilon \to 0$ and arbitrary  $\phi  \in
\Cvan(M)$, since, as a consequence of the bounded convergence theorem of Lebesgue for
Bochner integrals, 
\begin{displaymath}
\slim_{\epsilon \to 0} \epsilon^{-d} \pi(\chi_{Q_\epsilon})\phi= \slim_{\epsilon \to 0}
\int_{Q_1}  {\pi({\Phi^{-1}(\epsilon t)})\phi} \d t=\phi.
\end{displaymath}
Let $\phi\in \Cvan(M)_\infty$. From equation \eqref{eq:1} one deduces that $\chi_{Q_\epsilon}(t) A_j \pi(\Phi^{-1}(t)) \phi$ is Bochner integrable, so that using Lebesgue, and integrating by parts,         
we obtain 
\begin{align*}
\slim_{h\to0} h^{-1}(\pi(e^{ha_j})&-\1)\pi(\chi_{Q_\epsilon})\phi=\slim_{h\to0} h^{-1} \int_{Q_\epsilon} (\pi(e^{ha_j})-\1) \pi({\Phi^{-1}(t)})\phi\d t\\
&=\int_{Q_\epsilon} A_j \pi({\Phi^{-1}(t)})\phi\d t=\sum_{k=1}^d \int_{Q_\epsilon}\frac \gd{\gd t_k} \Big ( \pi({\Phi^{-1}(t)})\phi\Big ) d_{kj}(t) \d t\\
\intertext{}
&=\sum_{k=1}^d \int_{\hat Q_\epsilon}\Big[ d_{kj}(t)\pi({\Phi^{-1}(t)})
\phi\Big]_{(t_1,\dots,0,\dots,t_d)}^{(t_1,\dots,\epsilon ,\dots,t_d)} \d t_1
\wedge\hat{\dots}\wedge \d t_d \\
&-\sum_{k=1}^d\int_{Q_\epsilon} \pi({\Phi^{-1}(t)})\phi \frac \gd {\gd t_k} d_{kj}(t)
\d t,
\end{align*}
where $j=1,\dots,d$, and $d_{kj}(t)=\frac d{d h}s_k^j(0,t)$. The symbol $\hat\,$ indicates that integration over the variable $t_k$ is suppressed. As a consequence, 
$\pi(\chi_{Q_\epsilon}) \phi \in \D(A_i)$ for all   $\phi \in \Cvan(M)_\infty$.
Let us denote the difference on the right hand side of the last equality by $F_j(\phi)$. It  is defined for arbitrary $\phi \in \Cvan(M)$, and continuous in $\phi$ with respect to the strong topology in  $\Cvan(M)$. Since $\Cvan(M)_\infty$ is dense, 
and $\pi(\chi_{Q_\epsilon})$ is bounded, we obtain, by the closedness of the $A_i$, that $\pi(\chi_{Q_\epsilon}) \phi \in \D(A_i)$ for arbitrary $\phi \in \Cvan(M)$, and  
\begin{equation}
\label{eq:4}
A_j\pi(\chi_{Q_\epsilon}) \phi=F_j(\phi)
\end{equation}
for all   $\phi \in \Cvan(M)$. Assume now that $\phi \in \bigcap _{i=q,\dots,
d} \D(A_i)$. Then $\slim_{\epsilon \to 0} \epsilon^{-d}\pi(\chi_{Q_\epsilon})
\phi=\phi$ and we show that, similarly, 
\begin{equation*}
\slim_{\epsilon \to 0} \epsilon^{-d} A_j\pi(\chi_{Q_\epsilon})\phi=A_j \phi
\end{equation*}
for all    $j=q,\dots,d$. Writing $d_{kj}(t)
=\delta_{kj}+\sum_{|\alpha| \geq 1} c_\alpha^{kj} t^\alpha$  with complex coefficients $c_\alpha^{kj}$, we obtain with \eqref{eq:4}, by  substitution of variables,
\begin{align*}
 \epsilon^{-d} A_j\pi(\chi_{Q_\epsilon})\phi&=\sum_{k=1}^d \epsilon^{-1} \int_{\hat Q_1}
 \delta_{kj} \Big [\pi({\Phi^{-1}(t)})\phi \Big
 ]_{(t_1,\dots,0,\dots,t_d) \epsilon}^{(t_1,\dots,1,\dots,t_d) \epsilon} \d t_1
\wedge\hat{\dots}\wedge \d t_d\\
&+\sum_{k=1}^d  \int_{\hat Q_1}\sum_{l\not=k} \frac {\gd d_{kj}}{\gd
 t_l} (0)\Big [t_l \pi({\Phi^{-1}(t)})\phi \Big
 ]_{(t_1,\dots,0,\dots,t_d) \epsilon}^{(t_1,\dots,1,\dots,t_d) \epsilon} \d t_1
\wedge\hat{\dots}\wedge \d t_d\\
&+\sum_{k=1}^d \epsilon^{-1} \int_{\hat Q_1}\Big ( \frac {\gd d_{kj}}{\gd
 t_k} (0) \, \epsilon \, \pi(\Phi^{-1}(\epsilon t_1,\dots,\epsilon ,\dots \epsilon
 t_d))\phi-0\Big )\d t_1
\wedge\hat{\dots}\wedge \d t_d\\
&+\sum_{k=1}^d \epsilon^{-1} \int_{\hat Q_1}\sum_{|\alpha|\geq2} c_\alpha^{kj}\Big [t^\alpha \pi({\Phi^{-1}(t)})\phi \Big
 ]_{(t_1,\dots,0,\dots,t_d) \epsilon}^{(t_1,\dots,1,\dots,t_d) \epsilon} \d t_1
\wedge\hat{\dots}\wedge \d t_d\\
&-\sum_{k=1}^d  \int_{ Q_1}\frac {\gd d_{kj}}{\gd
 t_k} (\epsilon t)\pi(\Phi^{-1}(\epsilon t))\phi \d t.
\end{align*}
As $\epsilon$ goes to zero, the second and fourth summand on the right hand vanish, while the third and fifth cancel each other. Hence, only the first summand 
\begin{displaymath}
\epsilon^{-1}\int_{\hat Q_1} \pi(\e{\epsilon t_1a_1})\cdots [\pi(\e{\epsilon a_j})-\1]\cdots \pi(\e{\epsilon t_da_d})\phi \, dt_1
\wedge\hat{\dots}\wedge \d t_d
\end{displaymath}
remains. In the following, we prove that        
\begin{displaymath}
\slim _{\epsilon \to 0} \epsilon^{-1}\pi(\e{\epsilon t_1a_1})\cdots [\pi(\e{\epsilon a_j})-\1]\cdots
\pi(\e{\epsilon t_da_d})\phi=A_j \phi
\end{displaymath}
for all    $j=q,\dots,d$. For   $j=d$, the assertion is clear. Let, therefore, $j$ be equal $d-1$. Then,    
\begin{displaymath}
\epsilon^{-1} [\pi(\e{\epsilon a_{d-1}})-\1]\pi(\e{\epsilon t_da_{d}})\phi=\epsilon^{-1}
[\pi(\e{\epsilon a_{d-1}})-\1]\phi+ [\pi(\e{\epsilon a_{d-1}})-\1] \epsilon^{-1}[\pi(\e{\epsilon t_da_{d}})-\1]\phi.
\end{displaymath}
Now, since $\epsilon^{-1}[\pi(\e{\epsilon t_da_{d}})-\1]\phi$ converges strongly to $t_d A_d\phi$, and $\norm{\epsilon^{-1}[\pi(\e{\epsilon t_da_{d}})-\1]\phi}$
to    $\norm{t_d A_d\phi}$, in particular remaining bounded for $\epsilon \to 0$, we obtain
\begin{displaymath}
\norm{[\pi(\e{\epsilon a_{d-1}})-\1] \epsilon^{-1}[\pi(\e{\epsilon t_da_{d}})-\1]\phi} \to 0,
\end{displaymath}
and hence the assertion for   $j=d-1$. By iteration we get, for arbitrary $j=q,\dots,d$, 
\begin{gather*}
\epsilon^{-1}[\pi(\e{\epsilon a_j})-\1]\pi(\e{\epsilon t_{j+1}a_{j+1}})  \cdots \pi(\e{\epsilon t_da_d})\phi
\\=\frac {\pi(\e{\epsilon a_j})-\1} \epsilon \phi +[\pi(\e{\epsilon a_j})-\1]\sum_{i=j+1}^d \Big
(\prod_{m=j+1}^{i-1} \pi(\e{\epsilon t_{m}a_{m}})\Big )\frac {\pi(\e{\epsilon t_ia_i})-\1} \epsilon \phi,
\end{gather*}
and a repetition of the arguments above finally yields
the desired statement for $j=q,\dots,d$. Taking all together, we conclude for $\phi \in \bigcap _{i=q}^d \D(A_i)$ that
\begin{equation}
\label{eq:5}
\epsilon^{-d}(\pi(\chi_{Q_\epsilon})\phi,A_q\pi(\chi_{Q_\epsilon})\phi,\dots,A_d
\pi(\chi_{Q_\epsilon})\phi) \to (\phi,A_q\phi,\dots,A_d\phi)
\end{equation}
as $ \epsilon \to 0$ with respect to the strong product topology in $\Cvan(M)\times\dots\times
\Cvan(M)$. Now, for $\phi \in \Cvan(M)_\infty$,
\begin{align}
\label{eq:5a}
\begin{split}
(\pi(\chi_{Q_\epsilon})\phi,A_q\pi(\chi_{Q_\epsilon})&\phi,\dots,A_d
\pi(\chi_{Q_\epsilon})\phi)\\&=\int_{Q_\epsilon} (\pi({\Phi^{-1}(t)})\phi,A_q\pi({\Phi^{-1}(t)})\phi,
\dots, A_d\pi(\Phi^{-1}(t))\phi) \d t.
\end{split}
\end{align}
 If $\epsilon$ is sufficiently small, and 
$\phi \in \P$, then, by Proposition \ref{prop:2}, $$(\pi({\Phi^{-1}(t)})\phi,A_q\pi({\Phi^{-1}(t)})\phi,
\dots, A_d\pi(\Phi^{-1}(t))\phi) \in \overline
{\Gamma_{a_q,\dots,a_d|\P}},\qquad t \in Q_\epsilon,$$ so that the integral \eqref{eq:5a} also lies in $\overline
{\Gamma_{a_q,\dots,a_d|\P}}$. Since $\P$ is dense, it follows with  \eqref{eq:4}, $F_j$ being continuous, that 
\begin{displaymath}
(\pi(\chi_{Q_\epsilon})\phi,A_q\pi(\chi_{Q_\epsilon})\phi,\dots,A_d
\pi(\chi_{Q_\epsilon})\phi) \in \overline
{\Gamma_{a_q,\dots,a_d|\P}}
\end{displaymath}
for arbitrary $\phi \in \Cvan(M)$ and sufficiently small  $\epsilon$. With  \eqref{eq:5} we  finally obtain $ (\phi,A_q\phi,\dots,$ $A_d\phi)\in\overline{\Gamma_{a_q,\dots,a_d|\P}}$ for $\phi \in \bigcap_{i=q}^d \D(A_i)$. This proves the theorem.
\end{proof}

\section{Analytic elements of $(\pi,\Cvan(M))$}

Let $\pi$ be a continuous representation of a Lie group $G$ on a Banach space $\B$, and denote by $\B_\omega$ the space of all analytic elements in $\B$, i.\ e.\ the space of all
 $\phi \in \B$ for which $g\mapsto \pi(g)\phi$ is an analytic map from $G$ to $\B$. $\B_\omega$
is invariant under the  $G$-action   $\pi$ and the action  $d\pi$ of $\g$, norm-dense in  $\B$, and one has the inclusion 
$\B_\omega \subset \B_\infty$. Analytic elements of Banach representations were first studied by Harish--Chandra in \cite{harish-chandra},
and their importance is due to the fact that the closure of every $d\pi(\U(\g_\C))$-invariant subspace  of
$\B_\omega$ constitutes a $\pi(G)$-invariant subspace of $\B$. In addition note that every closed, $\pi(G)$-invariant subspace is also $\pi(C_c(G))$-invariant, where $C_c(G)$ denotes the space of continuous, complex valued functions on $G$ with compact support. Assume now that $G$ is a real linear algebraic group acting regularly on a smooth, real affine variety $M$, and consider the left regular representation of $G$ on $\Cvan(M)$ as introduced in the previous section.
The following theorem states that the elements of $\P= \C[M]\cdot e^{-r^2}$
are contained  in $\Cvan(M)_\omega$,  and is a consequence of the approximation argument given in Section \ref{sec:2.2}.

\begin{theorem}
\label{thm:2}
The elements of the $\g$-module $\P= \C[M]\cdot e^{-r^2}$ are analytic vectors of the left regular representation $\pi$ of  $G$
on  $\Cvan(M)$. 
\end{theorem}
\begin{proof}
Let $\phi=p \cdot e^{-r^2}\in \P$. Since $g \mapsto h^{-1}g$ is  an analytic isomorphism for arbitrary $h \in G$, it suffices to show that $g\mapsto \pi(g)\phi$ is analytic in a neighbourhood of the unit element, see e.\ g.\
\cite{harish-chandra}.  According to  \eqref{eq:6}, $\pi(g) \phi$ is given by the series
\begin{equation*}
 \pi(g) \phi=  \sum_{l=0}^\infty
\rho(g)p  (r^2-\rho(g)r^2)^l/l! \cdot e^{-r^2},
\end{equation*}
which converges absolutely in $\Cvan(M)$ provided that $g\in G$ is contained in a sufficiently small neighbourhood $U$ of  $e$. Since $G$ acts regularly on the affine variety $M$, $(g,m) \mapsto g  m$ is an analytic map from $G \times M$ to   $M$, and therefore  $(g,m)
\mapsto (\rho(g)p)(m)$  an analytic function, being also  polynomial in $m$. Thus,
\begin{equation}
\label{eq:7a}
(\rho(g)p) (m)=p(g^{-1}m)=\sum_{\Lambda,\beta} c_{\Lambda\beta}
\theta^\Lambda (g^{-1}) m^\beta,
\end{equation}
with constants $c_{\Lambda\beta}$ and multiindices $\Lambda,\beta$,
where
$\theta^\Lambda(g)=\theta_{11}^{\Lambda_{11}}(g)\dots
\theta_{nn}^{\Lambda_{nn}}(g)$, and $m^\beta=m_1^{\beta_1}\dots m_n^{\beta_n}$, the sum being finite. The coefficient functions  
$\theta_{ij}(g)=g_{ij}$, as well as their powers, are real analytic. In a neighbourhood of the unit element, they are given 
by the absolute convergent Cauchy product series 
\begin{equation}
\label{eq:7b}
\theta^\Lambda(g)=\sum_{\modulus{\gamma}\geq 0}^\infty b^\Lambda_{\gamma}
\Theta^{\gamma}(g)=\sum_{\gamma_1,\dots,\gamma_d\geq 0}^\infty b^\Lambda_{\gamma_1,\dots,\gamma_d}
\Theta_1^{\gamma_1}(g)\dots \Theta_d^{\gamma_d}(g),
\end{equation}
where the $\Theta_1,\dots \Theta_d$ are supposed to be coordinates near $e$ with 
$\Theta_i(e)=0$. By rearranging the sums we obtain 
\begin{equation*}
(\rho(g)p) (m)=\sum_{|\gamma|\geq 0}^\infty
p_{\gamma}(m)\Theta^{\gamma}(g^{-1}), 
\end{equation*}
where we put $p_{\gamma}(m)=\sum_{\Lambda,\beta} c_{\Lambda\beta}
b^\Lambda_{\gamma} m^\beta$. In a similar way, $(g,m)
\mapsto (r^2-\rho(g)r^2)^l(m)$ is a real analytic function on  $G\times M$ which is a polynomial expression in $m$. The coefficient functions have an expansion of the form $\theta_{ij}(g)=\delta_{ij}+\sum_{|\gamma|\geq 1}^\infty b^{ij}_{\gamma}
\Theta^{\gamma}(g)$, and we put  $\eta_{ij}(g)=\sum_{|\gamma|\geq 1}^\infty b^{ij}_{\gamma}
\Theta^{\gamma}(g)$. In this way we get
\begin{align*}
m^2&-(g^{-1}m)^2=\sum_{i=1}^n
(1-\theta_{ii}^2(g^{-1}))m_i^2-\sum_{i\not=j} ^n \theta^2_{ij}(g^{-1})
m_i^2-2\sum _{k=1}^n \sum_{i<j}^n
\theta_{ki}(g^{-1})\theta_{kj}(g^{-1})m_im_j\\
&=-\sum_{i,j=1}^n
\eta_{ij}^2(g^{-1})m_i^2-2\sum_{i\not=j}^n \eta_{ij}(g^{-1}) m_im_j-2\sum _{k=1}^n \sum_{i<j}^n
\eta_{ki}(g^{-1})\eta_{kj}(g^{-1})m_im_j,
\end{align*}
thus obtaining
\begin{equation}
\label{eq:7c}
\big (m^2-(g^{-1}m)^2\big )^l= \sum_{l\leq |\Lambda|\leq 2l,|\beta|=2l} d^l_{\Lambda
\beta} \eta^\Lambda(g^{-1}) m^\beta,
\end{equation}
where the $d^l_{\Lambda\beta}$ are real numbers, and the $\eta^\Lambda(g)$ are given by the Cauchy product series
\begin{equation}
  \label{eq:7d}
  \eta^\Lambda(g)=\sum_{\modulus{\delta}\geq \modulus{\Lambda}}^\infty b^\Lambda_\delta \Theta^\delta(g).
\end{equation}
 As a consequence, we can represent $(m^2-((g^{-1}m)^2)^l$ near the unit as an infinite sum in the variables $\Theta_1(g^{-1}),\dots,\Theta_d(g^{-1})$ in which only summands of order $\geq l$ do appear. Explicitly, one has     
\begin{equation*}
\big (m^2-(g^{-1}m)^2\big )^l=\sum_{|\delta|\geq l}^\infty q_{\delta}^l(m)\Theta^{\delta}(g^{-1}),
\end{equation*}
with  $q^l_\delta(m)=\sum_{l\leq|\Lambda|\leq \min(2l,\modulus{\delta}),|\beta|=2l} d^l_{\Lambda
\beta} b^\Lambda_\delta m^\beta$. For   
$|\delta| < l$ one has $q^l_\delta\equiv 0$. By choosing $U$ sufficiently small, we may assume that the coordinates $\Theta_1,
\dots,\Theta_d$ are defined on $U$. In the sequel, we shall show that, on $U$, $\pi(g)\phi$ has the expansion
\begin{equation}
\label{eq:7e}
\pi(g)\phi=\sum_{|\gamma|,|\delta|\geq
0}^\infty p_{\gamma}\Big(\sum_{l=0}^{|\gamma+\delta|}
q_{\delta}^l/l!\Big) \cdot e^{-r^2} \,\Theta^{\gamma+\delta}(g^{-1}).
\end{equation}
For this sake, we first compute 
\begin{align*}
\sum_{|\delta|\geq 1}^\infty |q^1_\delta(m) \Theta^\delta(g^{-1})|\leq\sum_{1\leq|\Lambda|\leq 2,|\beta|=2} |d^1_{\Lambda
\beta}| \Big (\sum_{|\delta|\geq \modulus{\Lambda}}^\infty |b^\Lambda_\delta
\Theta^\delta(g^{-1})| \Big )|m^\beta| \leq K_g m^2,
\end{align*}
where $K_g$ is a constant which depends only on $g$ and goes to zero for $g \to e$. But then, because of
$q^l_\delta(m)=\sum_{\delta^{(1)}+\dots+\delta^{(l)}=\delta}
q^1_{\delta^{(1)}}(m) \cdots q^1_{\delta^{(l)}}(m)$,
each summand of the series appearing in  \eqref{eq:7e} can be estimated according to
\begin{gather*}
\Big\|p_\gamma\cdot \Big(\sum_{l=0}^{|\gamma+\delta|}
q_\delta^l/l!\Big )\cdot e^{-r^2} \, \Theta^{\gamma+\delta}(g^{-1})\Big\|\\=
\sup_M\Big |p_\gamma(m)\Theta^{\gamma}(g^{-1})\Big (\sum_{l=0}^{|\gamma+\delta|} 
\sum_{\delta^{(1)}+\dots+\delta^{(l)}=\delta}
q_{\delta^{(1)}}^1(m)\Theta^{\delta^{(1)}}(g^{-1})\cdots
q_{\delta^{(l)}}^1(m)\Theta^{\delta^{(l)}}(g^{-1})/l!\Big ) e^{-m^2}
\Big|\\
\leq\sup_M\Big [|p_\gamma(m)\Theta^{\gamma}(g^{-1})|
\sum_{l=0}^{|\gamma+\delta|} (K_g m^2)^l/l!\,e^{-m^2}\Big
]\leq
\Big (\sum_{\Lambda,\beta}|c_{\Lambda\beta}b^\Lambda_\gamma\Theta^{\gamma}(g^{-1})|\sup_M|
m^\beta e^{-m^2/2}|\Big )\\
\cdot
\Big (\sum_{l=0}^{|\gamma+\delta|} K_g^l/l! \sup_M |m^{2l}
e^{-m^2/2}|\Big )\leq \Big (\sum_\Lambda
C_\Lambda|b^\Lambda_\gamma\Theta^{\gamma}(g^{-1})|\Big )\Big ( \sum_{l=0}^{|\gamma+\delta|} K_g^l
(2l)^l e^{-l}/l!\Big ),
\end{gather*}
where we put  $C_\Lambda=\sum_\beta |c_{\Lambda,\beta}| \sup_M |m^{\beta}
e^{-m^2/2}|$, and made use of \eqref{eq:6b}. The sum over $\Lambda$ appearing in the last line is a finite linear combination of  summands  of the convergent series  $\sum_{|\gamma|\geq 0}^\infty |b^\Lambda_\gamma\Theta^{\gamma}(g^{-1})|$; the summands of the series over $l$ behave  as 
$(2K_g)^l/\sqrt l$ for big $l$. Since $K_g<1/2$ for sufficiently small $U$, we deduce
 that the sum in  \eqref{eq:7e} converges absolutely, by the comparison criterium and the convergence of the geometric series. Therefore, by Riemann's theorem on the rearrangement of terms of absolutely convergent series, we can write this sum also in the form
\begin{displaymath}
\sum_{l=0}^\infty\sum_{|\gamma|,|\delta|\geq 0}^\infty  p_\gamma \, q_\delta^l/l!\cdot e ^{-r^2} \Theta^{\gamma+\delta}(g^{-1}).
\end{displaymath}
Now, by taking into account equations \eqref{eq:7a}-\eqref{eq:7d}, and the definition of the polynomials $p_\gamma$, $q^l_\delta$, one deduces that
\begin{gather*}
\Big\|\rho(g)p \, (r^2-\rho(g)r^2)^l\cdot e^{-r^2} -\sum_{0\leq|\gamma|,|\delta|\leq N}
p_\gamma\,
q_\delta^l\cdot e^{-r^2}\,  
\Theta^{\gamma+\delta}(g^{-1}) \Big \|\\\leq \sup_M \Big
| \sum_{\Lambda,\beta,\Lambda'\beta'} c_{\Lambda\beta} d_{\Lambda'\beta'}^l
 m^{\beta+\beta'} e^{-m^2} \Big ( \theta^{\Lambda}(g^{-1})\eta^{\Lambda'}
(g^{-1})-\sum_{0\leq|\gamma|,|\delta|\leq N}
b^\Lambda_\gamma b^{\Lambda'}_\delta \Theta^{\gamma+\delta}(g^{-1})\Big
)\Big |\\
\end{gather*}
goes to zero as $N \to \infty$, and we obtain   
\begin{equation*}
\rho(g)p^2  (r^2-\rho(g)r^2)^l\cdot e^{-r^2}=\sum_{|\gamma|,|\delta|\geq 0}^\infty  p_\gamma  \,q_\delta^l\cdot e ^{-r^2} 
\Theta^{\gamma+\delta}(g^{-1});
\end{equation*}
using equation  \eqref{eq:6} we get \eqref{eq:7e}.
Thus, for $\phi \in \P$, $\pi(g)\phi$ is analytic in a neighbourhood of the origin, as contended. 
\end{proof}

\section{Subquotients and reducing series of $\P$ and $\Cvan(M)$}

Assume that $G$ is a real reductive  group as defined in  \cite{wallach}. Let $K$ be a maximal compact subgroup of $G$ with Lie-Algebra $\k$, and
\begin{displaymath}
  \g=\k\oplus \p
\end{displaymath}
the corresponding Cartan decomposition of $\g$. If $\pi$ is a continuous representation of $G$ on a Banach space $\B$,  its restriction to $K$ defines a  continuous representation of $K$ on $\B$. Denote by $\hat K$ the set of all equivalence classes of finite dimensional irreducible representations of $K$. For each $\lambda \in \hat K$, let $\xi_\lambda$ be the character of $\lambda$, $d(\lambda)$ the degree of $\lambda$, and $\chi_\lambda=d(\lambda) \xi_\lambda$. We define  
\begin{equation*}
  P(\lambda)=\pi(\bar \chi_\lambda)=d(\lambda)\int_K \bar \xi_\lambda(k) \pi(k) \d k,
\end{equation*}
$\d k$ normalized Haar measure on $K$. $P(\lambda)$ is a continuous projection of $\B$ onto $\B(\lambda)=P(\lambda)\B$. Note that $\B(\lambda)$ is the isotypic $K$-submodule of $\B$ of type $\lambda$, i.\ e.\ it consists of all vectors, the linear span of whose $K$-orbit is finite dimensional, and splits into irreducible $K$-submodules of type $\lambda$, see \cite{warner}.

Let us resume the study of the left regular representation $(\pi,\Cvan(M))$, as introduced in section \ref{sec:3}, of a real linear algebraic group $G$ acting regularly on a smooth, real affine algebraic variety $M$, and of the $\g$-module $\P$. From now on, we will assume that $G$ is stable under transpose, so that $G$ is real reductive in the sense specified above. Note that $\R^n$ is endowed with a $K$-invariant scalar product in a natural way. By choosing appropriate coordinates, we can therefore assume that the distance function $r^2$ is invariant under the action of $K$. As in the previous sections, let $\rho$ be the left regular representation of $G$ on $\C[M]$. It is locally regular, so that $\P$ is contained in the subspace of differentiable, $K$-finite vectors of $(\pi,\Cvan(M))$. 
\begin{lemma}
\label{lem:2}
  $\P=\C[M] \cdot e^{-r^2}$ is a $(\g,K)$-module in the sense of Harish-Chandra and Lepowsky.
\end{lemma}
\begin{proof}
  Plainly, $\P$ is $\pi(K)$- and $d\pi(\g)$-invariant. Now, for general $\phi \in \Cvan(M)_\infty$, $X\in \g$, and $g\in G$ we have
  \begin{align}
\label{eq:9}
    \begin{split}
\pi(g)d\pi(X)\phi &=\lim_{h \to 0} h^{-1} \pi(g)[\pi(\e{hX})-1]\phi\\
&=\lim_{h\to 0} h^{-1} [\pi(\e{\Ad(g)hX})-1] \pi(g) \phi=d\pi(\Ad(g)X)\pi(g) \phi.    
    \end{split}  
  \end{align}
Let $\phi \in \P$, $E_\phi=\text{Span}\mklm{\pi(K) \phi}$. Since $(\rho,\C[M])$ is locally regular, $E_\phi$ is a finite dimensional subspace of $\P$.
Thus, $\P$ is a $(\g,K)$-module  in the sense of Harish-Chandra and Lepowsky, see e.\ g.\ \cite{wallach}.
\end{proof}
\begin{lemma}
\label{prop:3}
   Let $\P(\lambda)=\P \cap \Cvan(M)(\lambda)$. Then $\P=\sum_{\lambda \in \hat K} \P(\lambda)$, and  $\overline{\P(\lambda)}=\Cvan(M)(\lambda)$. 
\end{lemma}
\begin{proof}
 Since the linear span of the  $K$-orbit $\pi(K)[ p\cdot e^{-r^2}]$ of an arbitrary  element in $\P$ is finite dimensional,  we obtain the first assertion. The second assertion follows with Proposition 4.4.3.4. of \cite{warner}. 
\end{proof}
In what follows, we will call $\P$ also the \emph{underlying $(\g,K)$-module of $\Cvan(M)$}. Note that in general Hilbert representation theory, this role is occupied by the space of differentiable, $K$-finite vectors. $\P(\lambda)$ is the $\lambda$-isotypic component for $K$ of $\P$.
 \begin{remark}
  The fact that the sum $\sum_{\lambda\in \hat K} \P \cap \Cvan(M)(\lambda)$ is dense in $\Cvan(M)$ can also be seen by the following general argument. Fix $\phi\in \Cvan(M)$, and $\epsilon >0$. Since $\P$ is dense, there exists a polynomial $p=\sum _\alpha c_\alpha m^\alpha \in \C[M]$ such that $\|\phi-p\cdot e^{-r^2}\| \leq \epsilon/2$. For any finite subset $F$ in $\hat K$ we define $\bar \chi_F=\sum _{\lambda \in F} \bar \chi_\lambda$. As a consequence of the $K$-invariance of $r^2$, we have
  \begin{equation*}
    (P(\lambda)\, p\cdot e^{-r^2})(m)=\int _K \bar \chi_\lambda(k) p(k^{-1}m) \d k\cdot  e^{-r^2}=\sum _\alpha c_\alpha \int _K \bar \chi_\lambda(k) (k^{-1}m)^\alpha \d k\cdot  e^{-r^2},
  \end{equation*}
so that $P(\lambda) \,p\cdot e^{-r^2} \in \P \cap \Cvan(M)(\lambda)$. 
Now, according to a theorem of Harish-Chandra, for any differentiable vector $\psi \in \Cvan(M)_\infty$, the Fourier series
$ \sum _{\lambda\in \hat K} P(\lambda) \psi$
converges absolutely to $\psi$, see  Theorem 4.4.2.1 in \cite{warner}. For this reason, we can choose $F$ in such a way that $\big \|(\pi(\bar \chi_F)-1)p\cdot e^{-r^2}\big \| \leq \epsilon/2$. Hence,
\begin{displaymath}
  \big \|\pi(\bar \chi_F) p\cdot e^{-r^2}-\phi\big \| \leq \big \|(\pi(\bar \chi_F)-1)p\cdot e^{-r^2}\big \|+\big \|\phi-p\cdot e^{-r^2}\big \|\leq \epsilon.
\end{displaymath}
Since for any $F$, $\pi(\bar \chi_F)p\cdot e^{-r^2} \in \sum_{\lambda \in \hat K} \P \cap \Cvan(M)(\lambda)$, we obtain the desired assertion. 
\end{remark}
\begin{remark}
Let $T$ be a maximal torus of $K$.  The equivalence classes of finite dimensional irreducible representations of $K$ are in one-to-one correspondence with dominant integral and $T$-integral forms $\mu$  on the complexification $\t_\C$ of the Lie-Algebra of $T$. Let $\lambda_\mu$ be the class corresponding to $\mu$, and $(\pi_\mu,V^\mu)$ a representative in $\lambda_\mu$. Denote by  $M(\Lambda)$ the Verma module of weight $\Lambda\in\t_\C^\ast$, and $L(\Lambda)$ the unique, non zero irreducible subquotient of $M(\Lambda)$, see e.\ g.\ \cite{wallach}. Then the differential of $\pi_\mu$ is equivalent to the $\k_\C$-module $L(\mu)$. On the other hand, every finite dimensional simple $\k_\C$-module has a highest weight $\Lambda$, and is equivalent to $L(\Lambda)$. Thus, every $\lambda_\mu$-isotypic component for $K$ is an isotypic component for $\k_\C$ of type $L(\mu)$, and we get $\P=\sum_{\kappa \in \hat \k_\C} \P_\kappa$, where $\hat \k_\C$ denotes the set of all equivalence classes of finite dimensional simple $\k_\C$-modules. Hence $\P$ is a Harish-Chandra module for $\g_\C$ in the sense of \cite{dixmier}.  
\end{remark}
As already noted at the beginning of the previous section, Theorem \ref{thm:2} allows us to derive $\pi(G)$-invariant decompositions of  $\Cvan(M)$ from algebraic decompositions of $\P$ into $d\pi(\U(\g_\C))$-invariant subspaces. To begin with, 
note that the representation $\d\pi$ of  $\g$ on $\P$ is equivalent to a representation $\d\pi'$ of $\g$ on $\C[M]$ given by
\begin{equation}
\label{eq:9a}
  \d\pi'(X)=\d\rho(X) -[\d\rho(X) r^2], \qquad X \in \g.
\end{equation}
Define now in $\P$ the subspaces 
$$W_p=d\pi(\U(\g_\C))\text{ Span}\mklm{\pi(K) p\cdot e^{-r^2}}, \qquad  p\in \C[M].$$
\begin{lemma}
  $W_p$ is a $(\g,K)$-module in the sense of Harish-Chandra and Lepowsky.
\end{lemma}
\begin{proof}
Similar to the proof of Lemma \ref{lem:2}.
\end{proof}
Let $W_p(\lambda)$ be the $\lambda$-isotypic component for $K$ of $W_p$. Although, in general, $\dim W_p(\lambda)=\infty$, so that $W_p$ is not an admissible $(\g,K)$-module, any subquotient of $W_p$ turns out to be admissible. 
\begin{proposition}
  The subquotients $W_p/W_{d\pi'(X)p}$ are admissible $(\g,K)$-modules for any $p \in \C[M]$ and $X \in \U(\g_\C)$. 
\end{proposition}
\begin{proof}
 It suffices to consider the case $M=\R^n$, so that one has the grading $\C[\R^n]=\bigoplus_{l\geq 0} \C^l[\R^n]$ by the polynomial degree. Let $p \in \C[M]$, $X\in \U(\g_\C)$. Since $d\pi(X) \, p \cdot e^{-r^2}=[d\pi'(X) p] \cdot e^{-r^2}$, one has $W_{d\pi'(X)p} \subset W_{p}$ by \eqref{eq:9}. 
Denote by $Z$ the complement of  $W_{d\pi'(X)p} (\lambda)$ in $W_{p}(\lambda)$. Since, for $Y \in \g$, $d\pi'(Y): \C^l[\R^n] \rightarrow \C^l[\R^n] \oplus \C^{l+2}[\R^n]$, one deduces $W_p(\lambda)/ W_{d\pi'(X)p}(\lambda)\simeq Z \subset \bigoplus_{l=0}^{k-2} \C^l[\R^n] \cdot e^{-r^2}$, where $k$ is the degree of $d\pi'(X)p$. Thus, for arbitrary $\lambda \in \hat K$, $(W_p/W_{d\pi'(X)p})(\lambda)\simeq W_p(\lambda)/W_{d\pi'(X)p}(\lambda)$  is finite dimensional, and the assertion follows.
\end{proof}
By definition, $W_p\subset \Cvan(M)_\omega$ is invariant under $\U(\g_\C)$, so that its closure $\overline{W}_p$ constitutes a $G$-invariant subspace. Therefore, by restricting $\pi$ to $\overline{W}_p$, we obtain a Banach representation of $G$ on $\overline W_p$. Let us introduce now some terminology from general representation theory.
\begin{definition}
  Let $\pi$ be a Banach representation of a Lie group $G$ on a Banach space $\B$. Then by a \emph{reducing series} for $\pi$ we shall mean such a series for the integrated form of $\pi$, i.\ e.\ a nested family $\Bf$ of closed $\pi(\C_c(G))$-stable subspaces of $\B$ such that $\mklm{0}$ and $\B$ belong to $\Bf$, and with the property that if $\Af\subset \Bf$, then $\bigcap_{\A\in\Af} \A$ and $\overline{\bigcup_{\A\in \Af} \A}$ lie in $\Bf$. A \emph{Jordan-H\"{o}lder-series} for $\pi$ is a maximal reducing series for $\pi$.
\end{definition}
Let $\Bf$ be a reducing series for $\pi$. Two subspaces  $\B_1,\B_2\in \Bf$ with $\B_1\subset \B_2$ are said to be \emph{adjacent} if there is no $\tilde \B$ in $\Bf$ such that $\B_1 \subset \tilde \B\subset \B_2$, where the inclusions are understood to be strict.  
\begin{definition}
 Let $\pi$ be a Banach representation of $G$ on $\B$. Then a reducing series $\Bf$ for $\pi$ is said to be \emph{discrete} if, whenever $\B_1,\B_2\in\Bf$ and $\B_1 \subset  \B_2$, there are subspaces $\tilde \B_1,\tilde \B_2$ in $\Bf$ with $\B_1 \subset  \tilde \B_1 \subset\tilde \B_2 \subset  \B_2$ such that $\tilde \B_1$ and $\tilde \B_2$ are adjacent in $\Bf$.
\end{definition}
Now, as a consequence of Theorem \ref{thm:2}, we have the following proposition. Note that a similar statement also holds in case that  $G$ is not reductive.
\begin{proposition}
\label{prop:4}
 Let $p \in \C[M]$ and $\Bf=\mklm{ \overline{W}_p,\overline{W}_{d\pi'(X)p},  \overline{W}_{d\pi'(YX)p},\overline{W}_{d\pi'(ZYX)p},\dots}$, where $X,Y,Z, \dots$ are arbitrary elements in $\g$. Then $\Bf$ is a discrete reducing series for $(\pi,\overline{W}_p)$.
\end{proposition}
\begin{proof}
  To avoid undue notation, let as denote the elements of $\Bf$ by $\B_1,\B_2,\dots$. By construction and Theorem \ref{thm:2}, the elements of $\Bf$ are closed, $\pi(C_c(G))$-invariant subspaces of $\B_1$, and we have the inclusions $\B_1\supset \B_2 \supset \B_3\supset \dots$. Further, for $\Af \subset \Bf$, $\bigcap_{\B_i\in\Af} \B_i$ and $\overline{\bigcup_{\B_i \in \Af} \B_i}$ lie in $\Bf$. Thus $\Bf$ is a reducing series. Finally, by noting that every two $\B_{i}, \B_{i+1}$ are adjacent, we obtain the desired result. 
\end{proof}
\begin{lemma}
   Let $W_p(\lambda)=W_p \cap \Cvan(M)(\lambda)$. Then $W_p=\sum_{\lambda \in \hat K} W_p(\lambda)$ and  $\overline{W_p(\lambda)}=\overline{W_p}(\lambda)$. 
\end{lemma}
\begin{proof}
  Similar to the proof of Lemma \ref{prop:3}.
\end{proof}
\begin{lemma}
 Let $\Bf=\mklm{\B_1,\B_2,\dots}$ be a reducing series for $(\pi,\overline W_p)$ as in Proposition \emph{\ref{prop:4}}.
 Then the subquotients $\B_i/\B_j$, $i<j$, are $K$-finite Banach representations of $G$.
\end{lemma}
\begin{proof}
Let $(\pi,\B)$ be a Banach representation of a Lie group $G$, and $V$ a closed, $G$-invariant subspace of $\B$. Then
\begin{displaymath}
  P(\lambda) (\phi+ V)=\int \bar \chi_\lambda(k) (\pi(k) \phi +V) dk,
\end{displaymath}
so that $\phi +V \in (\B/V)(\lambda)$ implies $\phi \in \B(\lambda)$. One therefore gets an epimorphism from $(\B/V)(\lambda)$ onto $\B(\lambda)/V(\lambda)$ by setting $\phi +V \mapsto \phi +V(\lambda)$. Since $V(\lambda)\subset V$, this is clearly an isomorphism. Hence,  $(\B_i/\B_j)(\lambda)\simeq \B_i(\lambda)/\B_j(\lambda)$, where $i<j$. Let $Z$ be the algebraic complement of $\B_j(\lambda)$ in $\B_i(\lambda)$. 
Again, it suffices to consider the case $M=\R^n$. Then, by the preceeding lemma,
\begin{displaymath}
  \B_i(\lambda)/\B_j(\lambda)\simeq Z \subset \bigoplus_{l=0}^{k} \C^l[\R^n] \cdot e^{-r^2}
\end{displaymath}
for some $k$. Hence $\dim (\B_i/\B_j)(\lambda)<\infty$ for all $\lambda \in \hat K$, as contended.
\end{proof}
Assume now that $\h$ is a Cartan-subalgebra of $\g_\C$. 
Let $\Phi=\Phi(\g_\C,\h)$ be a root system with respect to $\h$,  $\g_\alpha=\mklm{X\in\g_\C: [Y,X]=\alpha(Y)X \text{ for all } Y\in \h}$ the root space corresponding to $\alpha\in \h^\ast$, and $\Phi^+$ a set of positive roots. As usual, we set
\begin{displaymath}
  \n^+=\bigoplus_{\alpha \in \Phi^+} \g_\alpha,\quad \n^-=\bigoplus_{\alpha \in \Phi^+} \g_{-\alpha}, \qquad \b=\h\oplus \n^+.
\end{displaymath}
Consider the underlying $(\g,K)$-module $\P$ of $\Cvan(M)$,  and define for $\lambda \in \h^\ast$ the corresponding weight space
\begin{displaymath}
  \P_\lambda=\mklm{\phi \in \P: d\pi(X) \phi =\lambda(X)\phi \text{ for all } X \in \h}.
\end{displaymath}
Note that $d\pi(\g_\alpha) \P_\lambda \subset \P_{\lambda+\alpha}$, and assume  that $\P_0$ is not empty. 
\begin{proposition}
  Let $p\cdot e^{-r^2} \in \P_\lambda$, and put $V_p=d\pi(\U(\g_\C)) p \cdot e^{-r^2}$, as well as $V_{\n^+,p}=\sum_{X\in \n^+} d\pi(\U(\g_\C)X)p\cdot e^{-r^2}$. Then $V_p/V_{\n^+,p}$ is a highest weight module of weight $\lambda$.
\end{proposition}
\begin{proof}
By construction, $d\pi(\n^+) p\cdot e^{-r^2}\in V_{\n^+,p}$, so that $p\cdot e^{-r^2}+V_{\n^+,p}$ is annihilated by $\n^+$ and therefore constitutes a highest weight vector. Since $V_p/V_{\n^+,p}$ is generated by $p\cdot e^{-r^2}+V_{\n^+,p}$ as a $\U(\g_\C)$-module, the assertion follows.   
\end{proof}
Let $M(\lambda)=\U(\g_\C) \otimes_{\U(\b)} \C_{\lambda}$ be the Verma module of weight $\lambda$ given by the above root system. Then, by general theory \cite{dixmier}, we have the following corollary.
\begin{corollary}
1. $V_p/V_{\n^+,p}=d\pi(\U(\n^-))(p\cdot e^{-r^2} +V_{\n^+,p})$. 

\noindent
2. $V_p/V_{\n^+,p}$ equals the algebraic direct sum of finite dimensional weight spaces.  

\noindent
3. $V_p/V_{\n^+,p}$ has a central character.

\noindent
4. There exists precisely one $\g$-homomorphism $\Psi$ from $M(\lambda)$ onto $V_p/V_{\n^+,p}$ such that $\Psi(1 \otimes 1)=p\cdot e^ {-r^2}+V_{\n^+,p}$.
\end{corollary}

Let us illustrate the above results by two examples.
\begin{example}
Consider  $G=\SL(2,\R)$, $K=\SO(2)$, acting on  $M=\R^2$ by matrices, and denote by $(\pi,\Cvan(\R^2))$ the corresponding regular representation of $G$ on the Banach space  $\Cvan(\R^2)$. As a basis for $\g$, take the matrices 
\begin{displaymath}
  a_1=\left (\begin{array}{cc} 0 & 1 \\ 0 & 0 \end{array} \right ), \qquad  a_2=\left (\begin{array}{cc} 0 & 0 \\ 1 & 0 \end{array} \right ), \qquad   a_3=\left (\begin{array}{cc}1 & 0 \\ 0 &-1 \end{array} \right ).
\end{displaymath}
Then $\k=\text{Span}_\R\mklm{a_1-a_2}$. Let $\Phi(\g_\C,\k_\C)=\mklm{0,\alpha}$ be a root system of $\g_\C$ with respect to the Cartan subalgebra $\k_\C$, where $\alpha \in \k_\C^\ast$ is defined by $\alpha (Y)=i(Y_{12}-Y_{21})$. Let $\Phi^+=\mklm{\alpha}$, so that $\n^\pm=\g_{\pm\alpha}=\text{Span}_\C\mklm{X^\pm}$, where $X^\pm=a_1+a_2\mp ia_3$. Put $H=a_1-a_2$.
 Consider in $\C[\R^2]$ the polynomials
 \begin{equation}
\label{eq:10}   
 r^2=m_1^2+m_2^2,\qquad  q_\pm=m_1\pm im_2,
 \end{equation}
as well as the identical polynomial $\1$. Let $\rho$ denote the regular representation of $G$ on  $\C[\R^2]$. Setting $\gd^\pm=\gd _{m_1}$ $\pm i\gd_{m_2}$, one computes 
$$d\pi(X^\pm)=\pm iq_\pm \gd ^\pm, \qquad d\pi(H)=m_1\gd_{m_2}-m_2\gd_{m_1}.
$$
Since identical expressions hold for $d\rho$, one has $d\rho(\n^+)q^k_+=0$, so that with respect to $\k_\C$, $(d\rho,\C^k[\R^2])$ constitutes a $(k+1)$-dimensional, irreducible $\g_\C$-module of highest weight $\alpha k/2$ with highest weight vector $q_+^k$. The polynomials $r^{2k}$ and $q_-^k$  are of weight zero, respectively $-\alpha  k/2$.  Define now in $\Cvan(M)$ the infinite dimensional, undecomposable subspaces
$$
W_1=d\pi(\U(\g_\C)) e^{-r^2}, \qquad W_{r^2}=d\pi(\U(\g_\C)) r^2 \cdot e^{-r^2},\qquad W_{q_\pm}=d\pi(\U(\g_\C))q_\pm\cdot e^{-r^2}.
$$ 
They are $(\g,K)$-modules in the sense of Harish-Chandra and Lepowsky, and by taking suitable subquotients one gets admissible $(\g,K)$-modules as well as and highest weight modules.
Note that  $X^+$, $X^-$ and  $H$ satisfy the  the commutation relations of the complexification of ${\mathfrak{sl}}(2,\C)$, so that, in particular, $[d\pi(X^+),d\pi(X^-)]=-4i \,d\pi(H)$. Since $d\pi(X^\pm):\C^l[\R^2]\cdot e^{-r^2}\rightarrow \C^l[\R^2]\cdot e^{-r^2}\oplus \C^{l+2}[\R^2]\cdot e^{-r^2}$ maps $\P_\lambda$ into $\P_{\lambda\pm\alpha}$, we obtain for the  $(\g,K)$-module $\P=\C[\R^2] \cdot e^{-r^2}$ the decomposition
\begin{equation*}
\P=W_\1 \oplus W_{r^2}\oplus W_{q_+}\oplus W_{q_-}.
\end{equation*}
Each of the summands is $d\pi(\U(\g_\C))$-invariant by construction, so that, by Theorem \ref{thm:2}, their closures must be $\pi(G)$-invariant. Since $\P$ is dense, we get for  $\Cvan(\R^2)$ a $G$-invariant decomposition according to
\begin{equation*}
\Cvan(\R^2)=\overline{\overline W_\1 \oplus \overline  W_{r^2}\oplus \overline W_{q_+}\oplus \overline W_{q_-}}
\end{equation*}
into undecomposable closed subspaces, each of them having discrete reducing series. The representation $(\pi,\Cvan(\R^2))$ itself is not of finite $K$-type and therefore not admissible.
\end{example}

\begin{example}
Consider again $G=\SL(2,\R)$, $K=\SO(2)$, but  now acting on $M=\g=\mathfrak{sl}(2,r)\simeq\R^3$ via the adjoint representation of $G$. Explicitly, we have the identification
\begin{displaymath}
\R^3 \ni m \mapsto X_m=\left ( \begin{array}{cc}m_2 & m_1+m_3 \\ m_1-m_3 &-m_2 \end{array}\right )\in \g.
\end{displaymath}
Similarly, we write $m_X$ for the point in $\R^3$ corresponding to $X\in\g$ under this identification.
Let $(\pi,\Cvan(\R^3))$ be the corresponding regular representation, and consider in $\C[\R^3]$ the polynomials
 \begin{equation*}  
 r^2=m_1^2+m_2^2+m_3^2, \qquad p_G=\det X_m=-m^2_1-m^2_2+m^2_3.
 \end{equation*}
 If not stated otherwise, we continue with the notation of the previous example.
 Since $ (d/dh) (\e{hY}m)_{|h=0}=m_{[Y,X_m]}$ for $Y \in \g$, one computes, using \eqref{eq:6c},
\begin{equation*}
d\pi(X^\pm)=\pm 2 i (m_3\gd_\mp + q_\mp \gd_{m_3}), \qquad d\pi(H)=2(m_2 \gd _{m_1} - m_1 \gd_{m_2}),
\end{equation*}
as well as identical expressions for $d\rho$. Note that $d\rho(X)p_G=0$ for all $X\in \g$, and $d\rho(\n^+)p_G^l q_-^k=0$, where $l,k$ are non negative integers. Hence, defining in $\C^k[\R^3]$ the finite dimensional subspaces
$$
V^{k,l}= d\rho(\U(\g_\C)) [ \,p_G^l\,q_-^{k-2l} ], \qquad  2l \leq k,
$$ 
we get highest weight $\g_\C$-modules of weight $(k-2l)\alpha$, see \cite{dixmier}. Since $\dim \C^k[\R^3]=\sum_{j=1}^{k+1} j$ and $\dim V^{k,l}=2(k-2l)+1$, one deduces that $\C^k[\R^3]$ decomposes according to 
\begin{displaymath}
\C^k[\R^3]=V^{k,0}\oplus V^{k,1}\oplus \dots \oplus V^{k,[k/2]}.
\end{displaymath}
If we therefore set $W_p=d\pi(\U(\g_\C)) p\cdot e^{-r^2}$, $p \in \C[\R^3]$, one gets for $\P$ the decomposition
\begin{equation*}
\P=\bigoplus _{l=0}^\infty  p_G^l\big [W_1\oplus W_{m_3}\oplus W_{r^2+m_3^2}\oplus \bigoplus_{k=1}^\infty 
\big (W_{q^k_+}\oplus W_{q^k_-}\big )\Big ]   
\end{equation*}
into $d\pi(\U(\g_\C))$-invariant subspaces of analytic vectors for $(\pi,\Cvan(\R^3))$, leading to  a corresponding decomposition of $\Cvan(\R^3)$ into $G$-invariant subspaces. Consider further  the smooth, affine $G$-varieties $N_\xi=\mklm{m \in \R^3: p_G=\xi}$, where $\xi\not=0$, and the underlying $(\g,K)$-module $\P[\xi]=\C[N_\xi] \cdot e^{-r^2}$ of $\Cvan(N_\xi)$.  One has  $\C[N_\xi]\simeq  \C[\R^3]/\I_{N_\xi}$, where $\I_{N_\xi}$ denotes the vanishing ideal of $N_\xi$.  Since $\I_{N_\xi}$ is generated by the polynomial $p_G-\xi$, one obtains in this case the decomposition 
\begin{equation*}
\P[\xi]=W_1\oplus W_{m_3}\oplus W_{r^2+m_3^2}\oplus \bigoplus_{k=1}^\infty 
\big (W_{q^k_+}\oplus W_{q^k_-}\big ).   
\end{equation*}
 Note that in this example neither $\P$,  nor $\P[\xi]$, are  longer finitely generated as  $(\g,K)$-modules.  
\end{example}

\providecommand{\bysame}{\leavevmode\hbox to3em{\hrulefill}\thinspace}
\providecommand{\MR}{\relax\ifhmode\unskip\space\fi MR }
\providecommand{\MRhref}[2]{%
  \href{http://www.ams.org/mathscinet-getitem?mr=#1}{#2}
}
\providecommand{\href}[2]{#2}

\end{document}